\documentclass[12pt]{amsart}
\usepackage{amscd,amsmath,amssymb}
\usepackage{graphics}
\usepackage{amsthm}

\DeclareMathOperator{\Mimo}{Mimo}
\DeclareMathOperator{\Mod}{Mod}
\DeclareMathOperator{\Ker}{Ker}
\DeclareMathOperator{\Hom}{Hom}

\DeclareMathOperator{\im}{Im}
\DeclareMathOperator{\End}{End}
\DeclareMathOperator{\Sl}{\mathcal{S}}

\DeclareMathOperator{\Rad}{Rad}
\DeclareMathOperator{\rep}{Rep}
\DeclareMathOperator{\lrep}{rep}
\DeclareMathOperator{\fin}{\mathfrak{s}}
\DeclareMathOperator{\lmod}{mod}

\input prepictex   \input pictex    \input postpictex

\def \arr#1#2{\arrow <2mm> [0.25,0.75] from #1 to #2}

\def\smallsq#1{\plot 0 0  0.#1 0  0.#1 0.#1  0 0.#1  0 0 /}
\def\s{\put{$\smallsq5$} at -.25 0 }
\def\ss{\multiput{$\smallsq5$} at -.25 -.25  -.25 .25 / }
\def\sad{\put{} at 0 -0.25 \circulararc -85 degrees from 0 .15 center at 0 0 \circulararc 85 degrees from 0 .15 center at 0 0 }
\def\happy{\put{} at 0 0.25 \circulararc 85 degrees from 0 -.15 center at 0 0 \circulararc -85 degrees from 0 -.15 center at 0 0 }
\def\longsad{\put{} at 0 -0.25 \ellipticalarc  axes ratio 2:1  180 degrees from .3 0 center at 0 0 }
\def\longhappy{\put{} at 0 0.25 \ellipticalarc  axes ratio 2:1  180 degrees from -.3 0 center at 0 0 }

\begin{document}

\theoremstyle{definition}
\newtheorem*{defn}{Definition}

\theoremstyle{plain}
\newtheorem{lemma}{Lemma}
\newtheorem{prop}[lemma]{Proposition}
\newtheorem*{claim}{Claim}
\newtheorem{thm}[lemma]{Theorem}
\newtheorem{cor}[lemma]{Corollary}
\newtheorem*{subclaim}{Subclaim}

\theoremstyle{remark}
\newtheorem*{prf}{Proof}
\newtheorem*{subproof}{Proof of Subclaim}

\title[The Auslander and Ringel-Tachikawa Theorem...]{The Auslander and Ringel-Tachikawa Theorem for Submodule Embeddings}\footnotetext{2000 MSC. Primary: 16G60, Secondary: 16G70, 20E15, 47A15} \author{Audrey Moore}\footnotetext{This research was partially supported by NSF DDEP grant \#0831369} \maketitle 

\begin{abstract} Auslander and Ringel-Tachikawa have shown that for an artinian ring $R$ of finite representation type, every $R$-module is the direct sum of finitely generated indecomposable $R$-modules.  In this paper, we will adapt this result to finite representation type full subcategories of the module category of an artinian ring which are closed under subobjects and direct sums and contain all projective modules.  In particular, the results in this paper hold for subspace representations of a poset, in case this subcategory is of finite representation type. \end{abstract}

\section{Introduction}

\vspace{5 mm}

Auslander, \cite[Corollary 4.8]{A}, Ringel and Tachikawa \cite[Corollary 4.4]{RT} showed that if $R$ is an artinian ring of finite representation type, then every $R$-module is the sum of finitely generated indecomposable $R$-modules.  In this paper, we are interested in a relative version of this result:  Let $\Mod R$ be a module category which may not be of finite type, and $\Sl \subseteq \Mod R$ a subcategory which has only finitely many finite length indecomposable objects up to isomorphy.  Under which circumstances is every object in $\Sl$ a direct sum of finite length indecomposable objects?

\vspace{5 mm}

\begin{thm}Let $R$ be an artinian ring, $\Sl$ a full subcategory of $\Mod R$ which is closed under direct sums and subobjects with $R_R\in \Sl$.  If there are only finitely many finite length indecomposable objects in $\Sl$, up to isomorphy, then every object in $\Sl$ is a direct sum of finitely generated indecomposable subobjects.  In particular, every indecomposable module in $\Sl$ has finite length.
\end{thm}

\vspace{5 mm}

The proof of Theorem 1 will be given in Section 4, where we will adapt a proof of the classical Auslander and Ringel-Tachikawa Theorem, which was given in lecture notes by W. Zimmermann \cite{Z2}.  As representations of posets are a motivating example, we will start with a slightly more restricted situation in Section 2, and use Auslander-Reiten theory to prove Proposition 2 below, a restatement of Theorem 1 in this situation in terms of Auslander-Reiten theory.  The benefit of this method of proof is that we obtain not only the result, but also information about the individual indecomposables and the almost split morphisms.

\vspace{5 mm}

\begin{prop}Let $R$ be a right Morita ring and $\Sl$ a full subcategory of $\Mod R$ which is closed under direct sums and summands.  If   $\Gamma$ is a finite component of the Auslander-Reiten quiver for $\fin = \Sl \cap \lmod R$ such that $\Gamma$ contains a projective generator, only endofinite modules, and each Auslander-Reiten sequence in $\Gamma$ is an Auslander-Reiten sequence in $\Sl$, then $\fin$ has Auslander-Reiten sequences in $\Sl$, $\Gamma$ is the Auslander-Reiten quiver for $\fin$ and every object in $\Sl$ is a direct sum of indecomposable modules in $\Gamma$.
\end{prop}

\vspace{5 mm}

Section 3 will focus on representations of posets with coefficients in a right Morita ring, illustrating the results from Section 2.  As an example we obtain the following application to invariant subspaces of linear operators.  This situation is of particular interest since the category of modules over the incidence algebra has infinite type, and hence does not satisfy the conditions of the classical version of the Auslander and Ringel-Tachikawa Theorem.

\vspace{5 mm}

\begin{prop}Let $k$ be any field, $V$ a possibly infinite dimensional $k$-vector space, $T:V\rightarrow V$ a $k$-linear operator which acts on $V$ with nilpotency index $2$, and let $V_1, V_2, V_3$ be subspaces of $V$ which are invariant under the action of $T$ and such that $V_1\subseteq (V_2\cap V_3)$.  Then $V$ has a direct sum decomposition $V=\bigoplus _i W_i$ into $T$-invariant vector spaces such that $V_j=\bigoplus _i V_j\cap W_i$ holds for $j=1,2,3$ and where each $W_i$ is isomorphic to exactly one of the 25 systems pictured in Section 3.
\end{prop}

\vspace{5 mm}

\emph{Related Results:}  For $P$ the one point poset and $\Lambda = \frac{\mathbf{Z}}{p^n}$ we are dealing with the problem posed by Birkhoff in 1934 to classify all subgroups of any finite abelian $p^n$-bounded group $B$, up to an automorphism of $B$.  For $n\leq 5$, Richman and Walker have shown that any subspace representation of $P$ is the direct sum of finitely generated indecomposable subspace representations of $P$ \cite{RW}.  This is particularly interesting in the case $n=4$ or $n=5$, where the module category of the incidence algebra is not of finite representation type.  Recently, C.M. Ringel \cite{R} has given a proof of the Auslander and Ringel-Tachikawa Theorem in the case where $R$ is an artin algebra which uses the Gabriel-Roiter measure and yields a method for splitting off direct summands.

\vspace{5 mm}

\emph{Acknowledgements:}  This paper is part of the author's Ph.D. thesis advised by Markus Schmidmeier at Florida Atlantic University.  The author would like to thank the representation theory group in Bielefeld, where this paper was completed, for their hospitality and the NSF for supporting this international exchange through DDEP grant \#0831369.

\section{Auslander-Reiten sequences and a result about indecomposables}
\vspace{5 mm}

In this section, let $R$ be a (two-sided) artinian ring, $\Mod R$ the category of right $R$-modules, $\Sl$ a full subcategory of $\Mod R$ which is closed under direct sums and summands with $R_R\in \Sl$, and $\fin = \Sl \cap \lmod R$, the full subcategory of $\Sl$ which has as objects the finitely generated modules in $\Sl$.  We will give sufficient conditions for $\fin$ to have Auslander-Reiten sequences with factorization property in $\Sl$ and deduce that if $\fin$ is of finite type and has endofinite indecomposables, then any object in $\Sl$ can be written as a direct sum of finite length indecomposable modules.  The terminology below regarding Auslander-Reiten sequences is adapted from \cite{AS}.

\vspace{5 mm}

Let $B,C\in \fin$.  A morphism $g:B\rightarrow C$ is \emph{right almost split in} $\Sl$ if $g$ is not a split epimorphism, and for any $h:X\rightarrow C$ which is not a split epimorphism, with $X\in \Sl$, the map $h$ factors through $g$.  Define left almost split in $\Sl$ dually.  Then we say $\fin$ \emph{has right almost split morphisms in} $\Sl$ if for any indecomposable object $C\in \fin$, there is an object $B\in \fin$ and an $f:B\rightarrow C$ which is right almost split in $\Sl$.  Define $\fin$ \emph{has left almost split morphisms in} $\Sl$ dually.  We say $\fin$ is \emph{functorially finite in} $\lmod R$ \emph{with respect to} $\Sl$ if $\fin$ is functorially finite in $\lmod R$, and the approximations have the factorization property for test objects in $\Sl$.

\vspace{5 mm}

A non-split exact sequence $0\rightarrow A\xrightarrow{f} B\xrightarrow{g} C\rightarrow 0$ with $A,B,C\in \fin$ is an \emph{Auslander-Reiten sequence in} $\Sl$ if $f$ is left almost split in $\Sl$ and $g$ is right almost split in $\Sl$.  We say $\fin$ \emph{has Auslander-Reiten sequences in} $\Sl$ if the following conditions hold: \begin{enumerate}
                                                                                                                                                                                                                                                                                                                 \item[$(i)$]The category $\fin$ has left and right almost split morphisms in $\Sl$.
			 \item[$(ii)$]For each indecomposable, non-projective $C\in \fin$, there is an Auslander-Reiten sequence in $\Sl$ ending in $C$ with objects in $\fin$.
			 \item[$(iii)$]For each indecomposable, non-injective $A\in \fin$, there is an Aus\-lan\-der-Reiten sequence in $\Sl$ beginning in $A$ with objects in $\fin$.
                                                                                                                                                                                                                                                                                                                \end{enumerate}

\vspace{5 mm}

\begin{prop}Let $\Sl$ be a full subcategory of $\Mod R$ which is closed under direct sums and summands such that $\fin = \Sl \cap \lmod R$ is functorially finite in $\lmod R$ with respect to $\Sl$.  If $\lmod R$ has Auslander-Reiten sequences in $\Mod R$, then $\fin$ has Auslander-Reiten sequences in $\Sl$.
\end{prop}

\vspace{5 mm}

\begin{prf}From Auslander and Smal\o, we know that since $\fin$ is functorially finite in $\lmod R$, $\fin$ has Auslander-Reiten sequences in $\fin$ \cite[section 4]{AS}.  We are now ready to verify requirements $(i)$ through $(iii)$.

$(i)$: Let $A \in \fin$ be indecomposable.  Since $A\in \lmod R$, there is a left almost split morphism $f:A\rightarrow B$ in $\Mod R$ with $B\in \lmod R$.  Let $b:B\rightarrow B'$ be a left approximation of $B$ in $\fin$.  Then the composition $b\circ f:A\rightarrow B'$ is left almost split in $\Sl$.  A dual construction yields right almost split morphisms in $\Sl$. 

$(ii)$: Let $C\in \fin$ be an indecomposable non-projective module.  Then from $(i)$, there is a right almost split morphism $g:B\rightarrow C$ in $\Sl$ with $B\in \fin$.  Since $B$ is a finite length module, we get a minimal version of this morphism by decomposing $B = B'\oplus B''$, with $g|_{B'}:B'\rightarrow C$ right minimal and $g|_{B''}=0$ \cite[I Theorem 2.2]{ARS}.  So $g|_{B'}$ is a minimal right almost split morphism, and we have an exact sequence $${\mathcal E} :0\rightarrow A\xrightarrow{f}B\xrightarrow{g|_{B'}}C\rightarrow 0$$  Consider the Auslander-Reiten sequence in $\fin$ ending in $C$: $0\rightarrow \tilde{A}\xrightarrow{\tilde{f}} \tilde{B} \xrightarrow{\tilde{g}} C\rightarrow 0$.  Since the minimal right almost split morphism in $\fin$ ending in $C$ is unique up to isomorphism, we get an isomorphism $h$ and a corresponding kernel map $k$ such that the diagram commutes \cite[V Proposition 1.4]{ARS}:
$$\beginpicture\linethickness1mm
  \setcoordinatesystem units <1.3cm,.9cm>
  \put{$B'$} at 0 0
  \put{$A$} at -1.5 0
  \put{$0$} at -3 0
  \put{$C$} at 1.5 0
  \put{$0$} at 3 0
  \put{$0$} at -3 -1.5
  \put{$\tilde{A}$} at -1.5 -1.5
  \put{$\tilde{B}$} at 0 -1.5
  \put{$C$} at 1.5 -1.5
  \put{$0$} at 3 -1.5
  \arr{-2.75 0}{-1.75 0}
  \arr{-1.25 0}{-.25 0}
  \arr{.25 0}{1.25 0}
  \arr{1.75 0}{2.75 0}
  \arr{-2.75 -1.5}{-1.75 -1.5}
  \arr{-1.25 -1.5}{-.25 -1.5}
  \arr{.25 -1.5}{1.25 -1.5}
  \arr{1.75 -1.5}{2.75 -1.5}
  \arr{-1.5 -.25}{-1.5 -1.2}
  \arr{0 -.25}{0 -1.2}
  \plot 1.45 -.25  1.45 -1.25 /
  \plot 1.55 -.25  1.55 -1.25 /
  \put{$f$} at -.75 .3
  \put{$g|_{B'}$} at .75 .3
  \put{$\tilde{f}$} at -.75 -1.2
  \put{$\tilde{g}$} at .75 -1.2
  \put{$k$} at -1.25 -.75
  \put{$h$} at .25 -.75
  \endpicture$$ 
Hence $k$ is an isomorphism, so we see that $A\in \fin$ and $\mathcal{E}$ is an Auslander-Reiten sequence in $\fin$.  Since $\fin$ has left almost split morphisms in $\Sl$, there is a left minimal almost split morphism $\hat{f}:A\rightarrow \hat{B}$ in $\Sl$ with $\hat{B}\in \fin$.  Then $\hat{f}$ is also a minimal left almost split morphism in $\fin$, and so $\hat{f}\cong f$ by uniqueness of minimal almost split morphisms \cite[V Proposition 1.5]{ARS}.  So $f$ is minimal right almost split in $\Sl$, and $\mathcal{E}$ is an Auslander-Reiten sequence in $\Sl$.

$(iii)$ Dual to the proof of $(ii)$. $_\Box$
\end{prf}

\vspace{5 mm}

Notice that the assumption $\lmod R$ has Auslander-Reiten sequences in $\Mod R$ must be stated in the previous proposition.  This assumption is satisfied whenever $R$ is an artin algebra \cite[Theorem 3.9]{A2}, but for $R$ an artinian ring this may not be the case \cite{Z}.

\vspace{5 mm}

We say an additive subcategory of $\Mod R$ is of \emph{finite (representation) type} if it has only finitely many finite length indecomposable objects, up to isomorphy.

\vspace{5 mm}

\begin{cor}  Let $\Sl$ be a full subcategory of $\Mod R$ which is closed under direct sums and summands with $R_R\in \Sl$, and such that $\fin = \Sl \cap \lmod R$ is of finite type and has Auslander-Reiten sequences in $\Sl$.  Then every non-zero object in $\Sl$ has a non-zero summand of finite length, and in particular each indecomposable in $\Sl$ has finite length. \end{cor}

\vspace{5 mm}

\begin{prf}  Let $X\in \Sl$ be non-zero.  Then there is a nonzero morphism $f:P\rightarrow X$ with $P$ an indecomposable projective. By assumption, $P\in \Sl$.  If $f$ is a split monomorphism, then $X$ has a direct summand isomorphic to $P$, and we are done.  So suppose $f$ is not a split monomorphism.  By assumption, $\fin$ has Auslander-Reiten sequences in $\Sl$, so there is a left almost split morphism starting at $P$, say $g_1:P\rightarrow B_1$.  Since $f$ is not a split monomorphism, it factors over $g_1$, so there is a map $f_1:B_1\rightarrow X$ such that $f=f_1\circ g_1$.  Also, since $B_1\in \fin$, it is finitely generated, and we can write it as a sum of indecomposables: $B_1=B_{1,1}\oplus B_{1,2}\oplus ...\oplus B_{1,n}$, and the maps $f_1$ and $g_1$ are given by the maps on the indecomposable summands, say $f_1=(f_{1,i})$, $g_1=(g_{1,i})$.  Since $f$ is non-zero, $f_1$ is nonzero, and there is some $i$ with $0\not= f_{1,i}:B_{1,i}\rightarrow X$, say $i_1$. Either $f_{1,i_1}$ is a split monomorphism and we are done, or we continue.  

Since $\fin$ is of finite type, there are only finitely many finitely generated indecomposables in $\fin$.  Let $m$ be the maximum length of these indecomposable $R$-modules.  By the Harada-Sai Lemma, this process must terminate after at most $2^m-1$ steps, so we get a split monomorphism $B_{j,k}\rightarrow X$ for some indecomposable $B_{j,k}\in \fin$, and we see that $X$ has a non-zero finite length direct summand.  $_\Box$ \end{prf}

\vspace{5 mm}

Before proving the main result of the section, it may be useful to recall that a \emph{Morita Ring} is an artinian ring such that the injective envelopes of simple right $R$-modules are finitely generated, and the notion of purity below.

\vspace{2 mm}

\begin{lemma}The following are equivalent for a short exact sequence $${\mathcal E} : 0\rightarrow A\xrightarrow{f} B\xrightarrow{g} C\rightarrow 0$$
\begin{enumerate}
\item For all $_RX$, $f\otimes 1_X$ is a monomorphism.
\item For all finitely presented $N_R$, $\Hom (N,g)$ is an epimorphism.
\end{enumerate}
\end{lemma}

If ${\mathcal E}$ satisfies these conditions, then ${\mathcal E}$ is called a \emph{pure exact sequence}, and $f(A)\subseteq B$ is a pure submodule.  So split exact sequences must be pure exact, and pure exact sequences are split exact if $C$ is finitely presented.

\vspace{5 mm}

We can now prove Proposition 2.

\vspace{5 mm}

\begin{prf}Since $\Gamma$ is finite and contains a projective generator, all indecomposables in $\fin$ occur in $\Gamma$, so $\fin$ is of finite type.  Since the almost split morphisms in $\Gamma$ are almost split in $\Sl$, $\fin$ has Auslander-Reiten sequences in $\Sl$ and $\Gamma$ is the Auslander-Reiten quiver for $\fin$.  Let $X\in \Sl$, and $\mathcal{Z}$ be the set of all (internal) direct sums of finite length indecomposable direct summands of $X$.  Then $\mathcal{Z}$ is non-empty and partially ordered by the direct sum relation.  Furthermore, by Azumaya's Theorem \cite[Theorem 12.6]{AF} each chain in $\mathcal{Z}$ will look like $$0\subseteq X_1=N_1\subseteq X_2= N_1\oplus N_2\subseteq X_3=N_1\oplus N_2\oplus N_3\subseteq \cdots$$ where the $N_i$'s are direct sums of distinct indecomposables in $\fin$.  Then $\cup _{i\in I} X_i = \bigoplus _{i\in I} N_i$ is an upper bound for the chain in $\mathcal{Z}$, so by Zorn's Lemma, $\mathcal{Z}$ has a maximal element, say $X'$, which is a pure submodule of $X$.  Since $X'$ is a direct sum of finitely many isomorphism types of modules of finite endolength, $X'$ itself has finite endolength, and hence is pure injective \cite[Chapter 4]{CB}.  Since $X'$ is a pure injective pure submodule, $X'$ is a summand of $X$.  So we have $X=X'\oplus X''$ with $X''$ having no non-zero finite length summand.  By Corollary 5, $X''$ must be $0$, and so $X=X'$ is a direct sum of finitely generated indecomposable submodules. $_\Box$
\end{prf}

\vspace{5 mm}

\section{Application to subspace representations of posets}

\vspace{5 mm}

In this section we will let $\Lambda$ be a right Morita ring.  Let $P$ be a poset, $P^*$ the quiver obtained from $P$ be adding a largest point, $\Lambda P^*$ the incidence algebra of $P$, which is a free $\Lambda$-module with basis $\{ (i,j):i\leq j$ in $P^*\}$ as discussed in \cite[Chapter 1]{A2S}.  Consider the full subcategory $\Sl = \rep_{\Lambda} P$, of $\Mod \Lambda P^*$ which has as objects the subspace representations of $P$.  First we will show that $\fin=\lrep _\Lambda P$ is functorially finite in $\lmod \Lambda P^*$ with respect to $\Sl$.  With this we are able to describe $\Sl = \rep _\Lambda P^*$ in all cases where $\fin$ is of finite type, as we will see in one example.
\vspace{5 mm}

Familiarity with the case of the one point poset will be useful in the proof of the next lemma, so for the moment, let $P=\bullet$ and $X\in \lmod \Lambda P^*$.  Then $X$ consists of a triple $(X_1, X_*; X_\alpha)$, where $X_1, X_*$ are $\Lambda$-modules and $X_\alpha :X_1\rightarrow X_*$.  Let $\overline{e}:\Ker X_\alpha \rightarrow I$ be the injective envelope of $\Ker X_\alpha$.  Using the injective factoring property, we can lift $\overline{e}$ to a map $e:X_1\rightarrow I$.  In \cite{RSC} Ringel and Schmidmeier introduced $\Mimo (X)=(X_1,X_*\oplus I_1; (X_\alpha, e)^T)$ and showed that the canonical map $\pi :\Mimo (X)\rightarrow X$, is a right approximation for $X$ in $\fin$.  By the construction, we see the factorization properties hold for test objects in $\Sl$.  So for the one point poset, $\lmod \Lambda P^*$ has right approximations in $\fin$ with respect to $\Sl$.  In the proof of the next lemma, we will generalize this construction for arbitrary posets.  Notice that the assumption that $R$ is a right Morita ring is needed so that the injective envelope of a finite length module is again of finite length.

\vspace{5 mm}

\begin{lemma} Let $P$ be a poset.  The subcategory $\fin = \lrep _\Lambda P$ is functorially finite in $\lmod \Lambda P^*$ with respect to $\Sl$.
\end{lemma}

\begin{prf}For ease of notation, we label the vertices of $P$ by natural numbers $1,2,...,n$ via some ordering of the points which extends from the partial ordering.  The notation $i<j$ will refer to the ordering of the poset.  Let $X\in \lmod \Lambda P^*$.  Then $X=(X_i,X_\alpha)$ where $X_i$ is a $\Lambda$-module at vertex $i$, and $X_\alpha$ is a map $X_\alpha:X_{s(\alpha)}\rightarrow X_{t(\alpha)}$.  For $i\leq j$ in $P$, and $p$ a path from $i$ to $j$, let $X_p$ be the composition of the corresponding maps $X_\alpha$.  Since the map $X_p$ is independent of the chosen path by the commutativity relations in $\lmod \Lambda P^*$, so we can refer to this morphism as $X_{ij}$.

With this set up, first we will show that $X\in \lmod \Lambda P^*$ has a left approximation in $\fin$ with respect to $\Sl$.  Let $L(X)=(L_i, L_\alpha )$ where $L_i=\im X_{i*}$ and each $L_\alpha:L_{s(\alpha)}\rightarrow L_{t(\alpha)}$ is the inclusion map.  Then $L(X)$ together with $l:X\rightarrow L(X)$ where $l=(X_{1*}, X_{2*}, \cdots , X_{n*}, 1_{X_*})$ fulfills the factorization property, and hence is a left approximation of $X$ in $\fin$ with respect to $\Sl$.

It remains to show that $X$ has a right approximation in $\fin$ with respect to $\Sl$.  For this, we generalize the construction of $\Mimo$ recalled above.  First, we define $\Mimo _k (X)$ for a vertex $k$ in $P$.  Let $(I_k$, $\overline{e_k})$ be the injective envelope of the kernel of $X_{k*}$, and let $e_k:X_k\rightarrow I_k$ be the extension of $\overline{e_k}$ via the injective factoring property.  Then $\Mimo _k (X) = (M_i, M_\alpha)$ where 
$$M_i=\begin{cases}
       X_i & i\leq k \\
       X_i\oplus I_k & i\nleq k
      \end{cases}$$
and for an arrow $\alpha :i\rightarrow j$ in $P$:
$$M_\alpha =\begin{cases}
             X_\alpha & j\leq k \\
             (X_\alpha ,e_k\circ X_{ik})^T & i\leq k, j\nleq k \\
             X_\alpha \oplus 1_{I_k} & i,j\nleq k
            \end{cases}$$

So for each $k$, $\Mimo _k (X)$ is a representation with a monomorphism starting at ver\-tex $k$ provided all arrows starting at successors of $k$ are monomor\-phisms in the rep\-re\-sen\-ta\-tion $X$.  Hence $$R(X)=\Mimo _1(\Mimo _2(\cdots (\Mimo_n (X))\cdots ))$$ has mo\-no\-mor\-phisms for each arrow, so $R(X)\in \fin$.  The module $R(X)$ together with the projection $r:R(X)\rightarrow X$ satisfies the factorization property required, since each map in the chain does: $$R(X)=\Mimo _1(\Mimo _2(\cdots (\Mimo _n(X))\cdots)) \rightarrow \cdots \rightarrow \Mimo _n(X) \rightarrow X$$ Hence $r:R(X)\rightarrow X$ is a right approximation of $X$ in $\fin$ with respect to $\Sl$. $_\Box$
\end{prf}

\vspace{5 mm}

As an example, let $k$ be any field, $\Lambda = \frac{k[T]}{T^2}$, and $P$ the poset with corresponding quiver:
$$
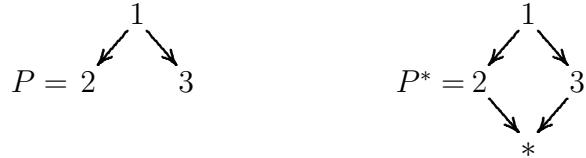
\beginpicture\linethickness1mm
  \setcoordinatesystem units <1.3cm,.9cm>
  \put{$P=$} at -3 0
  \put{$P^*=$} at 1 0
  \put{$1$} at -2 1
  \put{$2$} at -2.5 0
  \put{$3$} at -1.5 0
  \arr{-2.1 .75}{-2.4 .25}
  \arr{-1.9 .75}{-1.6 .25}
  \put{$1$} at 2 1
  \put{$2$} at 1.5 0
  \put{$3$} at 2.5 0
  \put{$*$} at 2 -1
  \arr{1.9 .75}{1.6 .25}
  \arr{2.1 .75}{2.4 .25}
  \arr{1.6 -.25}{1.9 -.75}
  \arr{2.4 -.25}{2.1 -.75}
  \endpicture$$

Then $X\in \lmod \Lambda P^*$ is given by the following information:

$$
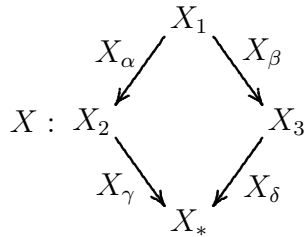
\beginpicture\linethickness1mm
  \setcoordinatesystem units <1.3cm,.9cm>
  \put{$X_3$} at -2.5 0
  \put{$X_1$} at -3.5 1.5
  \put{$X_2$} at -4.5 0
  \put{$X_*$} at -3.5 -1.5
  \arr{-3.75 1.25}{-4.25 .25}
  \arr{-3.25 1.25}{-2.75 .25}
  \arr{-4.25 -.25}{-3.75 -1.25}
  \arr{-2.75 -.25}{-3.25 -1.25}
  \put{$X:$} at -5.1 0
  \put{$X_\alpha$} at -4.25 1
  \put{$X_\beta$} at -2.75 1
  \put{$X_\gamma$} at -4.25 -1
  \put{$X_\delta$} at -2.75 -1
  \endpicture$$

To give a clear picture of the $\Mimo$ construction, below is the first step in computing $R(X)$:

$$\beginpicture\linethickness1mm
  \setcoordinatesystem units <1.3cm,.9cm>
  \put{$I_3\oplus X_2$} at 0 0
  \put{$X_1$} at 1 1.5  
  \put{$X_3$} at 2 0
  \put{$I_3 \oplus X_*$} at 1 -1.5
  \arr{.75 1.25}{.25 .25}
  \arr{1.25 1.25}{1.75 .25}
  \arr{.25 -.25}{.75 -1.25}
  \arr{1.75 -.25}{1.25 -1.25}
  \put{${{e_3\circ X_\beta}\choose{X_\alpha}}$} at 0 1.1
  \put{$X_\beta$} at 1.9 1
  \put{$1_{I_3}\oplus X_{\gamma}$} at 0 -1
  \put{${{e_3}\choose{X_\delta}}$} at 2 -1
  \put{$\Mimo_3(X):$} at -1.6 0
  \endpicture$$

And our right and left approximations $R(X)$ and $L(X)$ are:

$$\beginpicture\linethickness1mm
  \setcoordinatesystem units <1.3cm,.9cm>
  \put{$I_3\oplus X_2\oplus I_1$} at 0 0
  \put{$X_1$} at 1 1.5  
  \put{$X_3\oplus I_2\oplus I_1$} at 2 0
  \put{$I_3 \oplus X_4\oplus I_2\oplus I_1$} at 1 -1.5
  \arr{.75 1.25}{.25 .25}
  \arr{1.25 1.25}{1.75 .25}
  \arr{.25 -.25}{.75 -1.25}
  \arr{1.75 -.25}{1.25 -1.25}
  \put{\scriptsize{\stack <1pt> {{$e_3\circ X_\beta$}, $X_\alpha$, $e_1$}}} at 0 1.1
  \put{\bigg(} at -.4 1.1
  \put{\bigg)} at .4 1.1
  \put{\scriptsize{\stack <1pt> {$X_\beta$, $e_2\circ X_\alpha$, $e_1$}}} at 2 1.1
  \put{\bigg(} at 1.6 1.1
  \put{\bigg)} at 2.4 1.1
  \put{\footnotesize{${1_{i_3}\oplus{X_{\gamma}\choose e_2}\oplus 1_{I_1}}$}} at -.5 -1
  \put{\footnotesize{${{e_3}\choose{X_\delta}}\oplus 1_{I_2}\oplus 1_{I_1}$}} at 2.5 -1
  \put{$R(X):$} at -1.5 0
  \put{$\im X_\delta$} at -3.2 0
  \put{$\im (X_\gamma \circ X_\alpha)$} at -4 1.5
  \put{$\im X_\gamma$} at -4.8 0
  \put{$X_*$} at -4 -1.5
  \arr{-4.25 1.25}{-4.75 .25}
  \arr{-3.75 1.25}{-3.25 .25}
  \arr{-4.75 -.25}{-4.25 -1.25}
  \arr{-3.25 -.25}{-3.75 -1.25}
  \put{$L(X):$} at -6 0
  \put{$\iota$} at -4.75 1
  \put{$\iota$} at -3.25 1
  \put{$\iota$} at -4.75 -1
  \put{$\iota$} at -3.25 -1
  \endpicture$$

\vspace{5 mm}

Notice that $\Sl = \rep _\Lambda P$ is really the category of embeddings of embeddings of $\Lambda$-modules, and that the corresponding category of maps of maps, $\Mod \Lambda P^*$ is known to be of infinite type.  We compute a component $\Gamma$ of the Auslander-Reiten quiver in $\Sl$ using coverings, approximations (as above), and methods from \cite{RS}.  For a representation $X=(X_i)_{i\in P^*}$ we use the columns of boxes to denote the $\frac{k[T]}{T^2}$ module $X_*$ corresponding to the total space at $(*)$, $\frown$ the image of a generator of the submodule $X_2$, $\smile$ the image of a generator of $X_3$, and $\bullet$ the image of a generator of $X_1$.

$$\beginpicture\linethickness1mm
  \setcoordinatesystem units <.7cm,.7cm>
  \multiput{$\s$} at 2 2  3.75 4.25  5.75 6.25  8 12  8 8  9.75 5.75  12 4  14 2 /
  \multiput{$\ss$} at 0 4  0 8  0 12  1.75 6  2.25 6  4 12  4 8  4.25 4  4 0  6.25 6  5.75 2.25  6.25 1.75  7.75 4.25  8.25 3.75  8 0  10 2  10.25 6  12 0  12 8  12 12  14 6  14 10  14 14  16 4  16 8  16 12 /
  \multiput{$\happy$} at 4.25 4.25  1.75 5.75  6.25 6.25  0 7.75  4 8.25  12 7.75  16 8.25  14 10.25  0 12.25  4 11.75  8 12  16 11.75  8 .25  12 .25  10 2.25  14 2  0 3.75  12 4  16 3.75  14 5.75 /
  \multiput{$\sad$} at 6.25 2  8.25 4  2.25 6.25  1.75 5.75  9.75 5.75  0 8.25  4 7.75  8 8  16 7.75  0 11.75  4 12.25  12 11.75  16 12.25  14 14.25  8 .25  12 .25  10 2.25  14 2  0 3.75  12 4  16 3.75  14 5.75 /
  \multiput{\tiny$\bullet$} at 12 .25  10 1.75  14 2  0 3.75  8.25 3.5  16 3.75  1.75 5.75  6.25 5.75  4 7.75  4 11.75 /
  \multiput{$\longsad$} at 6 6.25  4 4.25 /
  \multiput{$\longhappy$} at 6 2  8 4  10 5.75  2 6.25 /
  \plot 0 14.5  0 12.65 /
  \plot 0 11.35  0 8.65 /
  \plot 0 7.35  0 4.65 /
  \plot 0 3.35  0 -.5 /
  \plot 16 14.5  16 12.65 /
  \plot 16 11.35  16 8.65 /
  \plot 16 7.35  16 4.65 /
  \plot 16 3.35  16 -.5 /
  \arr{.3 11.1}{1.7 6.9}
  \arr{.5 7.5}{1.3 6.7}
  \arr{.5 4.5}{1.3 5.3}
  \arr{.5 3.5}{1.5 2.5}

  \arr{2.3 6.9}{3.7 11.1}
  \arr{2.7 6.7}{3.5 7.5}
  \arr{2.7 5.3}{3.3 4.7}
  \arr{2.5 2.5}{3.5 3.5}
  \arr{2.5 1.5}{3.5 .5}

  \arr{4.3 11.1}{5.7 6.9}
  \arr{4.5 7.5}{5.3 6.7}
  \arr{4.7 4.7}{5.5 5.5}
  \arr{4.7 3.3}{5.3 2.7}
  \arr{4.5 .5}{5.5 1.5}

  \arr{6.3 6.9}{7.7 11.1}
  \arr{6.7 6.7}{7.5 7.5}
  \arr{6.7 5.3}{7.3 4.7}
  \arr{6.5 2.5}{7.5 3.5}
  \arr{6.7 1.3}{7.5 .5}

  \arr{8.3 11.1}{9.7 6.9}
  \arr{8.5 7.5}{9.5 6.5}
  \arr{8.5 4.5}{9.3 5.3}
  \arr{8.7 3.3}{9.5 2.5}
  \arr{8.5 .5}{9.5 1.5}

  \arr{10.3 6.9}{11.7 11.1}
  \arr{10.7 6.7}{11.5 7.5}
  \arr{10.7 5.3}{11.5 4.5}
  \arr{10.5 2.5}{11.5 3.5}
  \arr{10.5 1.5}{11.5 .5}

  \arr{12.5 12.5}{13.5 13.5}
  \arr{12.3 11.1}{13.7 6.9}
  \arr{12.5 8.5}{13.5 9.5}
  \arr{12.5 7.5}{13.5 6.5}
  \arr{12.5 4.5}{13.5 5.5}
  \arr{12.5 3.5}{13.5 2.5}
  \arr{12.5 .5}{13.5 1.5}

  \arr{14.5 13.5}{15.5 12.5}
  \arr{14.3 6.9}{15.7 11.1}
  \arr{14.5 9.5}{15.5 8.5}
  \arr{14.5 6.5}{15.5 7.5}
  \arr{14.5 5.5}{15.5 4.5}
  \arr{14.5 2.5}{15.5 3.5}

  \setdots<2pt>
  \plot .25 12  3.75 12 /
  \plot 4.25 12  7.75 12 /
  \plot 8.25 12  11.75 12 /
  \plot 0 2  1.75 2 /
  \plot 4.25 0  7.75 0 /
  \plot 14.25 2  16 2 /
  \plot .25 8  3.75 8 /
  \plot 4.25 8  7.75 8 /
  \plot 8.25 8  11.75 8 /

  \setshadegrid span <1mm>
  \vshade   0  2 12 <,z,,>
            2  2 12 <z,z,,>
            4  0 12 <z,z,,>
            8  0 12 <z,z,,>
           10  2 12 <z,z,,>
           12  0 12 <z,z,,>
           14  2 14 <z,z,,>
           16  2 12 /
  \endpicture$$ 

\vspace{2 mm}

For example, if 
$$\beginpicture\linethickness1mm
  \setcoordinatesystem units <1.3cm,.9cm>
  \put{$N=$} at 0 0
  \put{$\ss$} at .75 0
  \put{$\ss$} at 1.25 0
  \put{$\happy$} at .75 -.23
  \put{$\sad$} at .75 -.27
  \put{$\sad$} at 1.25 .23
  \put{\tiny$\bullet$} at .75 -.25
  \put{$\longhappy$} at 1 .28
  \put{$M=$} at 2.5 0
  \put{$\ss$} at 3.25 0
  \put{$\happy$} at 3.25 .25
  \put{$\sad$} at 3.25 .25
  \put{\tiny$\bullet$} at 3.25 .25
  \endpicture$$

Then $M$ stands for the representation with $M_1=M_2=M_3=M_*=\Lambda$ while $N$ stands for the representation with $N_* = \Lambda \oplus \Lambda$, $N_1 = k\oplus 0$, $N_2 = k\oplus \Lambda$, and $N_3 = \Lambda \cdot (1,1) \oplus k$ where $k$ denotes the socle of $\Lambda$.

\vspace{5 mm}

We can read off immediately that $\fin$ is of finite type since only 25 objects occur and $\Gamma$ contains a projective generator.  Since $\Lambda$ is an artin algebra, each indecomposable is endofinite, and $\lmod \Lambda P^*$ has Auslander-Reiten sequences in $\Mod \Lambda P^*$.  By Lemma 7 $\fin$ is functorially finite and hence has Auslander-Reiten sequences in $\Sl$ by proposition 4.  Thus the conditions of Proposition 2 are satisfied, and for this example each object in $\Sl$ is a direct sum of finite length indecomposable modules.  In other words, any indecomposable embedding of embeddings of $\frac{k[T]}{T^2}$-modules has finite length.  Also, as a consequence of Proposition 2, we obtain the following detailed information about the indecomposable summands of an arbitrary element of $\Sl$.

\vspace{5 mm}

\begin{cor}Let $X\in \rep _\Lambda P$ be an arbitrary representation.  Then $X$ has a direct sum decomposition into indecomposables such that each indecomposable summand $Y$ occurs in our list of 25 indecomposables and has the following properties:
\begin{itemize}
 \item $Y_*$ is a direct sum of at most two $\frac{k[T]}{T^2}$-modules.
 \item $Y$ is either isomorphic to $N$ or both $Y_2$ and $Y_3$ are cyclic as $k[T]$ modules.
 \item $Y$ is either isomorphic to $M$ or $\dim Y_1\leq 1$.
\end{itemize}
\end{cor}

\vspace{5 mm}

The above result can also be stated in terms of invariant subspaces of linear operators since $\frac{k[T]}{T^2}$-modules are pairs $(V,T)$ where $V$ is a vector space and $T:V\rightarrow V$ is a linear operator acting with nilpotency index 2.  Viewed in this way, submodules of $\frac{k[T]}{T^2}$-modules are invariant subspaces, and we see that Proposition 3 is just a translation of Corollary 7 into the language of invariant subspaces.

\vspace{5 mm}

\section{Proof of the main result}

\vspace{5 mm}

For this section, we require that $R$ is an artinian ring, $\Sl$ contains a projective generator and is closed under direct sums and subobjects such that $\fin =\lmod R\cap \Sl$ is of finite type.  Lemmas 9--14 leading up to a restatement of the Auslander and Ringel-Tachikawa Theorem have been copied or adapted from lecture notes of Zimmermann \cite{Z2} where the classical version of the Auslander and Ringel-Tachikawa Theorem is shown.  For the convenience of the reader, proofs are recalled.

\vspace{2 mm}

Let $M_1,\cdots M_n$ be a list of all indecomposables in $\fin$ up to isomorphy, and let $M=M_1\oplus \cdots \oplus M_n$.  Fix $X\in \Sl$, and let $H=\Hom_{R} (M,X)$.  We will denote $\Hom _R(Y,Z)$ by $(Y,Z)$.  Let $S=\End M_R$.  Then $M$ is a left $S$-module, denoted $_SM$, and $(M,X)_S$ is a right $S$-module via the structure $g\circ s(m)=g(s(m))$.

\vspace{5 mm}

\begin{lemma}  ${\displaystyle X=\sum _{h\in H} \im h}$
\end{lemma}

\begin{prf}  Let $X'$ be a finitely generated submodule of $X$. The category $\Sl$ is closed under submodules, so $X'\in \Sl$ and hence in $\fin$.  Since $\fin$ is of finite type, there are $k_1,\ldots k_n \in \mathbf{N}_0$ such that $X' \cong M_1^{k_1}\oplus \cdots \oplus M_n^{k_n}$.  So $X'\subseteq \displaystyle{\sum_{h\in H}\im h}$ and therefore $X=\displaystyle{\sum_{h\in H}\im h}$ $_\Box$
\end{prf}

\vspace{10 mm}

By this lemma, there is an epimorphism $p: M^{(H)} \rightarrow X$ given by $(m_h)_{h\in H} \mapsto \displaystyle{\sum _{h\in H} h(m_h)}$.  With $K=\Ker p$, we have a short exact sequence $${\mathcal E} :0\rightarrow K\rightarrow M^{(H)}\xrightarrow{p} X\rightarrow 0$$

\vspace{5 mm}

\begin{lemma}The sequence ${\mathcal E} :0\rightarrow K\rightarrow M^{(H)}\xrightarrow{p} X\rightarrow 0$ is pure exact.
\end{lemma}

\vspace{5 mm}

\begin{prf}Let $F\in \Sl$ be finitely presented.  We need to check that the sequence $$0\rightarrow (F,K) \rightarrow (F,M^{(H)}) \xrightarrow{(F,p)} (F,X) \rightarrow 0$$ is exact, in particular, that $(F,p)$ is an epimorphism.  Let $f\in (F,X)$.  Since $F$ is finitely generated, $F$ is a summand of $M^m$ for some $m\in \textbf{N}_0$.  So there is a split monomorphism $j:F\rightarrow M^m$ and a morphism $q:M^m \rightarrow F$ such that $1_F=q\circ j$.  Let $\iota _i:M\rightarrow M^m$ be the $i^{th}$ inclusion map.  Then we have the following diagram.
$$
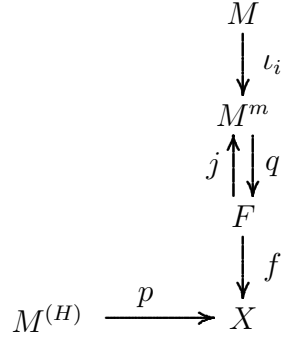
\beginpicture\linethickness1mm
  \setcoordinatesystem units <1.3cm,.9cm>
  \put{$M$} at 0 0
  \put{$M^m$} at 0 -1.5
  \put{$F$} at 0 -3
  \put{$X$} at 0 -4.5
  \put{$M^{(H)}$} at -2 -4.5
  \arr{0 -.3}{0 -1.2}
  \put{$\iota _i$} at .3 -.75
  \arr{.1 -1.8}{.1 -2.7}
  \arr{-.1 -2.7}{-.1 -1.8}
  \put{$q$} at .3 -2.25
  \put{$j$} at -.3 -2.25
  \arr{0 -3.3}{0 -4.2}
  \put{$f$} at .3 -3.75
  \arr{-1.4 -4.5}{-.3 -4.5}
  \put{$p$} at -1 -4.2
  \endpicture$$

Let $f_i=f\circ q\circ \iota _i :M\rightarrow X$.  Then $f_i \in H$ is the restriction of $p$ to the summand of $M^{(H)}$ with index $f_i$, by definition of $p$.  Denote by $\iota _{f_i}$ the canonical inclusion $\iota _{f_i}:M\rightarrow M^{(H)}$ into the $f_i$th component in $M^{(H)}$.  Let $g:M^m\rightarrow M^{(H)}$ be the map induced by the $\iota_{f_i}$ so that $g\circ \iota _i = \iota _{f_i}$. Finally,let $h=g\circ j:F\rightarrow M^{(H)}$.  Then we have the following commutative diagram.
$$
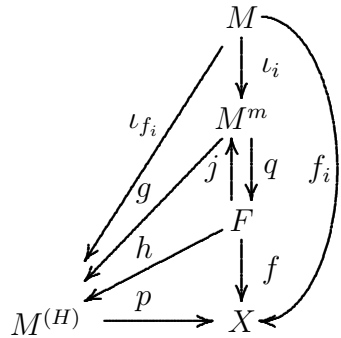
\beginpicture\linethickness1mm
  \setcoordinatesystem units <1.3cm,.9cm>
  \put{$M$} at 0 0
  \put{$M^m$} at 0 -1.5
  \put{$F$} at 0 -3
  \put{$X$} at 0 -4.5
  \put{$M^{(H)}$} at -2 -4.5
  \arr{0 -.3}{0 -1.2}
  \put{$\iota _i$} at .3 -.75
  \arr{.1 -1.8}{.1 -2.7}
  \arr{-.1 -2.7}{-.1 -1.8}
  \put{$q$} at .3 -2.25
  \put{$j$} at -.3 -2.25
  \arr{0 -3.3}{0 -4.2}
  \put{$f$} at .3 -3.75
  \arr{-1.4 -4.5}{-.3 -4.5}
  \put{$p$} at -1 -4.2
  \arr{-.2 -.45}{-1.6 -3.6}
  \put{$\iota _{f_i}$} at -1 -1.6
  \arr{-.2 -1.8}{-1.6 -3.9}
  \put{$g$} at -1 -2.6
  \arr{-.2 -3.15}{-1.6 -4.2}
  \put{$h$} at -1 -3.4
  \ellipticalarc axes ratio 1:2 180 degrees from .2 -4.5 center at .2 -2.25
  \arr{.2 -4.5}{.19 -4.5}
  \put{$f_i$} at .8 -2.25
  \endpicture$$

Since $f_i = fq\iota _i = p\iota _{f_i}$, we have $pg\iota _i = p\iota_{f_i}=f_i=fq\iota _i$, for all $i$.  So $pg=fq$, and $ph=pgj=fqj=f$. So for $f\in (F,K)$, there is an $h\in (F,M^{(H)})$ with $ph=f$; and $(F,p)$ is onto. $_\Box$
\end{prf}

\vspace{10 mm}

\begin{lemma} The sequence obtained by applying $(M,-)$ to ${\mathcal E}$, $$0\rightarrow (M,K) \rightarrow (M,M^{(H)}) \rightarrow (M,X) \rightarrow 0,$$  is a pure exact sequence of $S$-modules.
\end{lemma}

\vspace{5 mm}

\begin{prf}Since $0\rightarrow K\rightarrow M^{(H)}\rightarrow X\rightarrow 0$ is a pure exact sequence of $R$-modules, $$0\rightarrow (M,K) \rightarrow (M,M^{(H)}) \rightarrow (M,X) \rightarrow 0,$$ is an exact sequence of $S$-modules by Lemma 10 since $M$ is finitely presented.  It remains to show that the sequence is pure exact.

Let $N$ be a finitely generated right $S$-module.  Then there is an epimorphism $\pi :S^n \rightarrow N$ for some $n$.  Since $-\otimes _S M: \Mod S\rightarrow \Mod R$ is a right exact functor, we get an $R$ epimorphism $\pi \otimes 1:S^n\otimes _S M \rightarrow N \otimes _S M$.  Since tensor products commute with direct sums, $S^n\otimes _SM=M^n$ is a finitely generated $R$-module.  Hence $N\otimes _SM_R=\im (\pi\otimes 1)$ is also a finitely generated $R$-module.

Now we are ready to apply $(N,-)$ to our sequence and check that it gives a short exact sequence.

$$0\rightarrow (N_S,(M_R,K_R)_S)\rightarrow (N_S,(M_R,M_R^{(H)})_S)\rightarrow (N_S, (M_R,X_R)_S)$$

Notice that for a right $S$-module $N$, an $S$-$R$ bimodule $M$, and a right $R$-module Y, we have the adjoint isomorphism which is natural in $N$ and $Y$: 
\begin{eqnarray*}(N\otimes_S M,Y) \rightarrow (N,\Hom _R(M,Y)) \\
\phi \mapsto (n\mapsto (m\mapsto \phi (n\otimes m)) \end{eqnarray*}
So we obtain the commutative diagram:

$$\beginpicture\linethickness1mm
  \setcoordinatesystem units <1.1cm,.9cm>
  \put{$(N_S,(M_R,K_R)_S)$} at -4.5 0
  \put{$0$} at -7 0
  \put{$(N_S,(M_R,M^{(H)}_R)_S)$} at -.75 0
  \put{$(N_S,(M_R,X_R)_S)$} at 3 0
  \put{$0$} at -7 -2
  \put{$(N_S\otimes M_R,K_R)$} at -4.5 -2
  \put{$(N_S\otimes M_R,M^{(H)}_R)$} at -.75 -2
  \put{$(N_S\otimes M_R,X_R)$} at 3 -2
  \arr{-6.75 0}{-6 0}
  \arr{-3 0}{-2.5 0}
  \arr{1 0}{1.5 0}
  \arr{-6.75 -2}{-6 -2}
  \arr{-3 -2}{-2.25 -2}
  \arr{.75 -2}{1.5 -2}
  \arr{-4.5 -1.5}{-4.5 -.5}
  \arr{-.75 -1.5}{-.75 -.5}
  \arr{3 -1.5}{3 -.5}
  \put{$\cong$} at -4.25 -1
  \put{$\cong$} at -.5 -1
  \put{$\cong$} at 3.25 -1
  \endpicture$$

But $0\rightarrow K\rightarrow M^{(H)} \rightarrow X\rightarrow 0$ is pure exact and $N_S\otimes M_{R}$ is finitely generated, so the bottom sequence is short exact.  Hence $(N_S,(M_R,M^{(H)}_R)_S)\rightarrow (N_S,(M_R,X_R)_S)$ is an epimorphism by the communtativity of the diagram; and so the lemma is proved. $_\Box$
\end{prf}

\vspace{10 mm}

\begin{lemma} The evaluation map $\lambda :(M,X)\otimes _S M\rightarrow X_{R}$ where $\phi \otimes m \mapsto \phi (m)$ is an isomorphism of $R$-modules.
\end{lemma}

\vspace{5 mm}

\begin{prf}Let $x\in X_{R}$.  Then $x=\sum \phi _i(m_i)$ for some $\phi _i \in (M,X)$ and $m_i\in M$ by Lemma 9.  It remains to show that $\lambda $ is injective.  Let $\displaystyle{\sum_{i=1}^{m} \phi _i\otimes x_i} \in \Ker \lambda$, with $\phi _i\in (M,X)$, $x_i\in M$.  Then $\lambda \big(\displaystyle{\sum_{i=1}^{m} \phi _i\otimes x_i}\big) = \displaystyle{\sum_{i=1}^{m}\phi_i(x_i)} =0$.   Also, since $R_R\in \Sl$, Lemma 9 implies we can write $1=\displaystyle{\sum_{j=1}^n \alpha _j(m_j)}$ for some $\alpha _j\in (M,R)$, $m_j\in M$.  Let $s_{ij}:M\rightarrow M$ be the $R$-homomorphism given by $m\mapsto x_i\alpha _j(m)$.  Then 
\begin{eqnarray} \displaystyle{\sum_{i=1}^m \phi _i\otimes x_i} & = & \displaystyle{\sum_{i=1}^m \phi _i \otimes x_i\cdot 1}\\
&=& \displaystyle{\sum_{i=1}^m \phi _i\otimes x_i\sum_{j=1}^n\alpha _j(m_j)}\\
&=& \displaystyle{\sum_{i,j}\phi _i\otimes s_{ij}(m_j)} \\
&=& \displaystyle{\sum_{i,j}\phi _is_{ij}\otimes m_j}\\
&=& \displaystyle{\sum_j\left(\sum_i\phi _i s_{ij}\right)\otimes m_j}\\ 
&=& \displaystyle{\sum_j0\otimes m_j}\\
&=&0 \end{eqnarray}
Note that equality $(4)$ holds since $s_{ij}\in S$, and we can see that $(6)$ holds since for $x\in M$, $\displaystyle{\sum _i \phi _is_{ij}(x)=\sum_i\phi _i(x_i \alpha_j(x))}$ but $\alpha _j(x)$ is a scalar and $\displaystyle{\sum_{i=1}^{m}\phi_i(x_i)} =0$ by the kernel assumption.  So $\Ker \lambda =0$; and we see that $\lambda $ is a monomorphism and hence an isomorphism, as required. $_\Box$
\end{prf}

\vspace{10 mm}

\begin{lemma} $(M,M^{(H)})\cong S^{(H)}$ as right $S$-modules.
\end{lemma}

\vspace{5 mm}

\begin{prf} Since $S = \End(M)$, $S^{(H)}=\bigoplus _H (M,M)$.  Also, $(M,M^{(H)}) = (M, \bigoplus _H M)$.  But $M$ is fin\-it\-ely gen\-er\-at\-ed, so $\bigoplus _H (M,M) = (M,\bigoplus _H M)$. $_\Box$
\end{prf}

\vspace{10 mm}

In particular, since $S^{(H)}$ is a projective $S$-module; $(M,M^{(H)})\cong S^{(H)}$ is a projective $S$-module.  To prove the main result, however, we would like $(M,X)$ to be projective.  Since $S$ is the endomorphism ring of a finite length module, it is semiprimary; and so we will use the following general lemma.

\vspace{5 mm}

\begin{lemma} Let $S$ be semiprimary.  Let $$ {\mathcal E} :0\rightarrow A\xrightarrow{f} B\xrightarrow{g} C\rightarrow 0$$ be a pure exact sequence of $S$-modules with $B$ projective.  Then $C$ is also projective.
\end{lemma}

\vspace{5 mm}

\begin{prf} Since $S$ is semiprimary, every $S$-module has a projective cover, and each projective in $S$ is the direct sum of indecomposable projective $S$-modules. \cite[Theo\-rem 27.11]{AF}  Let $\pi :P\rightarrow C$ be a pro\-ject\-ive cover of $C$, $\displaystyle{P=\bigoplus _{i\in I} P_i}$ for some set $I$ and indecomposable projective $S$-modules $P_i$. We want to show that $L=\Ker (\pi )$ is zero, so pick $x\in L$.  Since $x\in L\subseteq P=\displaystyle{\bigoplus_{i\in I}P_i}$, there is a finitely generated projective summand $P'$ of $P$, with $x\in P'$.  Since $\pi:P\rightarrow C$ is an epimorphism and $g:B\rightarrow C$ with $B$ projective, the projective factoring property yields a map $t:B\rightarrow P$ such that $\pi t=g$.  So for the kernel map $t'$, we have the following commutative diagram with exact rows.
$$\beginpicture\linethickness1mm
  \setcoordinatesystem units <1.3cm,.9cm>
  \put{$B$} at 0 0
  \put{$A$} at -1.5 0
  \put{$0$} at -3 0
  \put{$C$} at 1.5 0
  \put{$0$} at 3 0
  \put{$0$} at -3 -1.5
  \put{$L$} at -1.5 -1.5
  \put{$P$} at 0 -1.5
  \put{$C$} at 1.5 -1.5
  \put{$0$} at 3 -1.5
  \arr{-2.75 0}{-1.75 0}
  \arr{-1.25 0}{-.25 0}
  \arr{.25 0}{1.25 0}
  \arr{1.75 0}{2.75 0}
  \arr{-2.75 -1.5}{-1.75 -1.5}
  \arr{-1.25 -1.5}{-.25 -1.5}
  \arr{.25 -1.5}{1.25 -1.5}
  \arr{1.75 -1.5}{2.75 -1.5}
  \arr{-1.5 -.25}{-1.5 -1.25}
  \arr{0 -.25}{0 -1.25}
  \plot 1.55 -.25  1.55 -1.25 /
  \plot 1.45 -.25  1.45 -1.25 /
  \put{$f$} at -.75 .25
  \put{$g$} at .75 .25
  \put{$\iota$} at -.75 -1.25
  \put{$\pi$} at .75 -1.25
  \put{$t'$} at -1.25 -.75
  \put{$t$} at .25 -.75
  \endpicture$$ 
Let $j:P'\rightarrow P$ be the inclusion map.  Since $x\in P'$, $xS$ is a submodule of $P'$,so let $\iota ':xS\rightarrow P'$ be this inclusion map.  The we have 
$$\beginpicture\linethickness1mm
  \setcoordinatesystem units <1.3cm,.9cm>
  \put{$B$} at 0 0
  \put{$A$} at -1.5 0
  \put{$0$} at -3 0
  \put{$C$} at 1.5 0
  \put{$0$} at 3 0
  \put{$0$} at -3 -1.5
  \put{$L$} at -1.5 -1.5
  \put{$P$} at 0 -1.5
  \put{$C$} at 1.5 -1.5
  \put{$0$} at 3 -1.5
  \arr{-2.75 0}{-1.75 0}
  \arr{-1.25 0}{-.25 0}
  \arr{.25 0}{1.25 0}
  \arr{1.75 0}{2.75 0}
  \arr{-2.75 -1.5}{-1.75 -1.5}
  \arr{-1.25 -1.5}{-.25 -1.5}
  \arr{.25 -1.5}{1.25 -1.5}
  \arr{1.75 -1.5}{2.75 -1.5}
  \arr{-1.5 -.25}{-1.5 -1.25}
  \arr{0 -.25}{0 -1.25}
  \plot 1.55 -.25  1.55 -1.25 /
  \plot 1.45 -.25  1.45 -1.25 /
  \put{$f$} at -.75 .25
  \put{$g$} at .75 .25
  \put{$\iota$} at -.75 -1.25
  \put{$\pi$} at .75 -1.25
  \put{$t'$} at -1.25 -.75
  \put{$t$} at .25 -.75
  \put{$0$} at -3 -3
  \put{$xS$} at -1.5 -3
  \put{$P'$} at 0 -3
  \put{$\frac{P'}{xS}$} at 1.5 -3
  \put{$0$} at 3 -3
  \arr{-2.75 -3}{-1.75 -3}
  \arr{-1.25 -3}{-.25 -3}
  \arr{.25 -3}{1.25 -3}
  \arr{1.75 -3}{2.75 -3}
  \arr{0 -2.75}{0 -1.75}
  \put{$j$} at .25 -2.25
  \put{$\iota '$} at -.75 -2.75
  \put{$\pi '$} at .75 -2.75
  \endpicture$$ 
Since $x\in L$, $xS\subseteq L$, and there is an inclusion map $j':xS\rightarrow L$, which makes the bottom left square commute and hence induces a cokernel map $j'':\frac{P'}{xS}\rightarrow C$.  Also, $P'$ is a summand of $P$, and $j$ is a split monomorphism, so there is a $q$ such that $1_{P'}=qj$.
$$\beginpicture\linethickness1mm
  \setcoordinatesystem units <1.3cm,.9cm>
  \put{$B$} at 0 0
  \put{$A$} at -1.5 0
  \put{$0$} at -3 0
  \put{$C$} at 1.5 0
  \put{$0$} at 3 0
  \put{$0$} at -3 -1.5
  \put{$L$} at -1.5 -1.5
  \put{$P$} at 0 -1.5
  \put{$C$} at 1.5 -1.5
  \put{$0$} at 3 -1.5
  \arr{-2.75 0}{-1.75 0}
  \arr{-1.25 0}{-.25 0}
  \arr{.25 0}{1.25 0}
  \arr{1.75 0}{2.75 0}
  \arr{-2.75 -1.5}{-1.75 -1.5}
  \arr{-1.25 -1.5}{-.25 -1.5}
  \arr{.25 -1.5}{1.25 -1.5}
  \arr{1.75 -1.5}{2.75 -1.5}
  \arr{-1.5 -.25}{-1.5 -1.25}
  \arr{0 -.25}{0 -1.25}
  \plot 1.55 -.25  1.55 -1.25 /
  \plot 1.45 -.25  1.45 -1.25 /
  \put{$f$} at -.75 .25
  \put{$g$} at .75 .25
  \put{$\iota$} at -.75 -1.25
  \put{$\pi$} at .75 -1.25
  \put{$t'$} at -1.25 -.75
  \put{$t$} at .25 -.75
  \put{$0$} at -3 -3
  \put{$xS$} at -1.5 -3
  \put{$P'$} at 0 -3
  \put{$\frac{P'}{xS}$} at 1.5 -3
  \put{$0$} at 3 -3
  \arr{-2.75 -3}{-1.75 -3}
  \arr{-1.25 -3}{-.25 -3}
  \arr{.25 -3}{1.25 -3}
  \arr{1.75 -3}{2.75 -3}
  \arr{.12 -2.75}{.12 -1.75}
  \put{$j$} at .37 -2.25
  \put{$\iota '$} at -.75 -2.75
  \put{$\pi '$} at .75 -2.75
  \arr{-.12 -1.75}{-.12 -2.75}
  \put{$q$} at -.37 -2.25
  \arr{-1.5 -2.75}{-1.5 -1.75}
  \arr{1.5 -2.7}{1.5 -1.75}
  \put{$j'$} at -1.25 -2.25
  \put{$j''$} at 1.75 -2.25
  \endpicture$$ 
But $\frac{P'}{xS}$ has projective resolution $P_1\rightarrow P'\xrightarrow{\pi '} \frac{P'}{xS}$, where $P_1$ is finitely generated since it is the projective cover of $\Ker \pi '=xS$.  Since $\frac{P'}{xS}$ is also finitely generated, it is finitely presented.  Also, ${\mathcal E}$ is pure exact, so for $j''\in \Hom (\frac{P'}{xR},C)$ there is a $\tau \in \Hom (\frac{P'}{xR}, B)$ such that $j''=g\tau$.  Since $\pi (j-t\tau \pi ')=j''\pi '-j''\pi '=0$ by commutativity, $j-t\tau \pi '\in \Ker \pi = \im \iota$, and there is a $\sigma :P'\rightarrow L$ so that $j-t\tau \pi '=\iota \sigma$.  But $\iota \sigma \iota '=(j-t\tau \pi ')\iota '=j\iota '-t\tau \pi '\iota '=j\iota '=\iota j'$ by commutativity and exactness of the third row, and $\iota \sigma \iota '=\iota j'$ implies $\sigma \iota '=j'$ since $\iota$ is injective.

Now let $h=\iota \sigma q\in \End P$.  Since $\pi h=\pi \iota \sigma q=0$, $\im h\subseteq \Ker \pi$.  But $\pi$ is a projective cover, so $\Ker \pi$ is small in $P$, and hence $\im h$ is also small in $P$.  So $h\in \Rad (\End (P))$.  Since $S$ is semiprimary; $\Rad (\End (P))$ is nilpotent; so $h$ is nilpotent, and $1-h$ is an automorphism of P.  Finally, $(1-h)(j\iota' (x))=j\iota '(x)-\iota \sigma qj\iota '(x)=j\iota '(x)-\iota \sigma \iota '(x)=0$.  Since $1-h$ is an automorphism; $j\iota '(x)=0$.  Since $j\iota '$ is a monomorphism; $x=0$, so $L=0$.  Hence $P\cong C$ and we see $C$ is projective. $_\Box$
\end{prf}

\vspace{10 mm}

Since $M$ is a finitely generated right $R$ module in $\Sl$, and $S=\End(M_R)$ we have the functor $$(-\otimes _S M_R): {\mathbf{P}}_S\rightarrow {\mathbf{S}}_M$$ where $\mathbf{P}_S$ is the category of projective right $S$-modules and $\mathbf{S}_M$ is the category of direct summands of direct sums of copies of $M$ \cite[Lemma 29.4]{AF}.  We can now complete the proof of Theorem 1:

\vspace{5 mm}

\begin{prf} Let $X\in \Sl$.  Let $M_1,\cdots M_n$ be a list of all indecomposables in $\mathcal{S}$ up to isomorphy, define $M=M_1\oplus \cdots \oplus M_n$, and $S=\End (M)$.  Then the short exact sequence $$0\rightarrow (M,K)\rightarrow (M,M^{(H)})\rightarrow (M,X)\rightarrow 0$$ is pure exact by Lemma 11.  By Lemma 13, the middle term is projective, and Lemma 14 yields that $(M,X)_S$ is projective.  Applying the functor $(-\otimes _S M_R)$ to $(M,X)_S \in {\mathbf{P}}_S$, we see $(M,X)_S\otimes _S M_R\in {\mathbf{S}}_M$.  By Lemma 12, $(M,X)_S\otimes _S M_R \cong X$, so $X\in {\mathbf{S}}_M$ and hence is a direct summand of $M^{(J)}$ for some set $J$.  By definition of $M$, $X$ is a direct summand of $M_1^{(J)}\oplus \cdots \oplus M_n^{(J)}$.  So by Azumaya's theorem; $X=M_1^{(J_1)}\oplus \cdots \oplus M_n^{(J_n)}$ where $J_1, \cdots J_n \subseteq J$.  $_\Box$
\end{prf}

\end{document}